\documentclass{amsart}
\usepackage{mathrsfs}
\usepackage{amsfonts}
\usepackage{amssymb}
\usepackage[papersize={6.9in,9.8in},textwidth=14cm,textheight=21.7cm,centering]{geometry}

\usepackage[colorlinks=true,citecolor=blue,linkcolor=blue]{hyperref}

\usepackage{hyperref}

\allowdisplaybreaks

\newtheorem*{theorem*}{Theorem}

\newtheorem{theorem}{Theorem}

\theoremstyle{remark}

\begin{document}

\title{On a problem of Bart Goddard}

\thanks{The work is supported by N.S.F. (No. 11171265) of P.R. China.}

\author{Ping Xi}

\address{School of Mathematics and Statistics, Xi'an Jiaotong University, Xi'an 710049, P. R. China}

\email{xprime@163.com}

\begin{abstract}We find identities for two series proposed by Bart Goddard in terms of elementary functions.\end{abstract}

\maketitle

\section{Introduction and the main result}

During the 2009 Western Number Theory held at Asilomar, Bart Goddard proposed a problem to ask the expressions of the series
\[A_k(y):=\sum_{n\geqslant1}(-1)^{n+1}\binom{2n+1}{2k}\frac{y^{2n-1}}{(2n+1)!},\]
\[B_k(y):=\sum_{n\geqslant1}(-1)^{n+1}\binom{2n+1}{2k+1}\frac{y^{2n-2}}{(2n+1)!}\]for any fixed natural number $k$. See \cite{M} [Problem 009:14] for details.

In this short note, we give the desired expressions for $A_k(y)$ and $B_k(y)$. To this end, we consider the following series
\begin{align}\label{eq:1}S_k(y)=\sum_{n\geqslant1}(-1)^{n+1}\binom{2n+1}{k}\frac{y^{2n+1}}{(2n+1)!}.\end{align}
Clearly, $A_k(y)=y^{-2}S_{2k}(y)$ and $B_k(y)=y^{-3}S_{2k+1}(y)$ for $y\neq0$. As usual, we regard the series above as the formal power series, and we don't consider their convergence in $y$.

In general, we can prove that
\begin{theorem}\label{thm:1}We have
\[S_0(y)=y-\sin y,\]\[S_1(y)=y-y\cos y,\]and for $k\geqslant2,$ we have
\[S_k(y)=\frac{y^k}{k!}\times\begin{cases}(-1)^{(k+2)/2}\sin y,\ \ \ & k\text{ even},\\
(-1)^{(k+1)/2}\cos y,\ \ \ & k\text{ odd}.\end{cases}\]
\end{theorem}

The proof is an easy application of the method of generating functions, which will be outlined in the next section.

\bigskip
\section{Proof of the theorem}
Now we put
\begin{align}\label{eq:2}S(x,y)=\sum_{k\geqslant0}S_k(y)x^k.\end{align}

In view of (\ref{eq:1}) and changing the order of summation, we have
\begin{align*}S(x,y)&=\sum_{n\geqslant1}(-1)^{n+1}\frac{y^{2n+1}}{(2n+1)!}\sum_{k\geqslant0}\binom{2n+1}{k}x^k=-\sum_{n\geqslant1}(-1)^{n}\frac{(xy+y)^{2n+1}}{(2n+1)!}.\end{align*}
Hence we have
\[S(x,y)=xy+y-\sin(xy+y)=xy+y-\sin xy\cos y-\cos xy\sin y.\]
Using the Taylor expansions for $\sin x$ and $\cos x$, we get
\begin{align}S(x,y)&=\sum_{k\geqslant0}\beta_k(y)\frac{x^k}{k!},\label{eq:3}\end{align}
where $\beta_0(y)=y-\sin y,\beta_1(y)=y-y\cos y$, and
\[\beta_k=\begin{cases}(-1)^{(k+2)/2}{y^k}\sin y,\ \ \ & k\geqslant2\text{ even},\\
(-1)^{(k+1)/2}{y^k}\cos y,\ \ \ & k\geqslant2\text{ odd}.\end{cases}\]

By comparing the coefficients in (\ref{eq:2}) and (\ref{eq:3}), we arrive at
\begin{align*}S_k(y)=\frac{\beta_k(y)}{k!},\end{align*}which completes the proof of the theorem.

\bigskip
\bigskip

\bibliographystyle{plainnat}

\bigskip

\bigskip

\end{document}